\documentclass[letterpaper, 10pt, journal, twoside]{IEEEtran}

\IEEEoverridecommandlockouts

\usepackage[latin1]{inputenc}
\usepackage{graphicx} 
\usepackage{float}		
\usepackage{cite}
\usepackage{color,soul}
\usepackage{calc}

\usepackage{amsmath,amsfonts,amssymb,amsthm}
\usepackage{mathtools}
\usepackage[american,smartlabels]{circuitikz}
\ctikzset{bipoles/length=0.7cm}
\ctikzset{bipoles/thickness=0.7}
\ctikzset{/tikz/circuitikz/nodes width=0.1}
\tikzstyle{every node}=[font=\small]
\tikzstyle{every path}=[line width=0.8pt,line cap=round,line join=round]

\newcommand{\Ex}{\mathcal{E}}
\definecolor{blu}{rgb}{0,0,1}

\definecolor{red}{rgb}{1,0,0}

\begin{document}

\title{Data-Driven LQR Control Design}

\author{Gustavo~R.~Gon\c{c}alves~da~Silva, 
        Alexandre~S.~Bazanella,
				Charles Lorenzini and
        Luc\'{i}ola~Campestrini
\thanks{The authors are with the Department of Automation and Energy -- Federal University of Rio Grande do Sul. Av. Osvaldo Aranha, 103 -- CEP: 90035-190 -- Porto Alegre, RS -- Brazil}
\thanks{E-mails \{gustavo.rgs,~bazanella,~charles.lorenzini,~luciola\}@ufrgs.br.}
\thanks{This work was supported by the National Council for Scientific and Technological Development -- CNPq/BR and by the Brazilian Federal Agency for Support and Evaluation of Graduate Education -- CAPES/BR. }}

\IEEEaftertitletext{This work has been submitted to the IEEE for possible publication. Copyright may be transferred without notice, after which this version may no longer be accessible.}

\maketitle

\begin{abstract}
This paper presents  a data-driven solution to the discrete-time infinite horizon LQR problem.
The state feedback gain is computed directly  from a batch of input and state data collected from the plant.
Simulation examples illustrate the convergence of the proposed solution to the optimal LQR gain as the number of Markov parameters tends to infinity. Experiments in an uninterruptible power supply are presented, which demonstrate the
practical applicability of the design methodology.
\end{abstract}

%
%

\section{INTRODUCTION}

\IEEEPARstart{T}{he} 
Linear Quadratic Regulator (LQR) design is a classical control problem whose analysis and solution can be found
in most  textbooks on control theory. It consists in computing the state feedback gain that optimizes
a quadratic cost function of the plant's state and input. This computation of the gain
requires the solution of a Riccati equation and
is given as a function of the plant's state-space model. Closed-form (also called \textit{batch-form}) solutions to the Riccati
equation have also been provided \cite{lewis1981,furuta1993}, and these are given as a function
of the plant's Markov parameters. 
Whether applying the classical approach of explicitly solving the Riccati equation, or using a plant's state space
description to calculate the Markov parameters and then feed them into the closed-form solution,
this is a {\em model-based} design approach. That is, it is a design approach that is based on the knowledge
of a good enough explicit model of the plant and on the use of this model
in the control design following the {\em certainty equivalence} principle.  


Data-Driven (DD) optimal control design methods have also been developed based on these closed-form solutions
of the Riccati equation  \cite{furuta1995, LQRmarkov, LQRestimado}. Although these DD design methods start from
the LQR/LQG problem formulation, they do not calculate the state feedback control gain; instead, they
directly estimate from data the optimal control input at each time instant. As such, they can not be said to
solve the LQR problem in its classical formulation, and can mainly be cast within a predictive control framework.

Motivated by applications in which a state feedback is to be designed but a good enough model is not available
and is of no interest \textit{per se}, we present in this paper a DD approach to the solution of the LQR design.
Otherwise stated, we provide a DD solution for the computation of the optimal state feedback gain. 
In a DD control design, the controller structure is defined \textit{a priori} and the controller's parameters are tuned
with the use of a large batch of data, usually after such data are acquired.
In most of the DD control literature, the controller structure consists of output feedback with a predefined transfer function
with parameters to be tuned -- such as a PID controller, for instance.
The design itself is based, in most methods, on the Model Reference approach \cite{Bazanella&Campestrini&Eckhard:2012}. 
A control design based on the DD approach leads naturally to the automation of the design process,
thus being extremely convenient for auto-tuning and self-tuning, and (all things being equal)
also tends to outperform  model-based designs, as shown in \cite{Campestrini:Eckhard:Bazanella:Gevers:2017}. 

The infinite-horizon LQR problem fits this formulation perfectly: one has a fixed controller structure (the state feedback gain) with a few parameters to tune, and the controller must minimize a given quadratic performance criterion.
The plant's model is just an intermediate step in the design and often has no interest in itself, the controller being
the final and only objective.  Thus, in this paper we  provide a method to compute the LQR state feedback gain from data without
the intermediate step of identifying a model of the system. 

The paper is organized as follows. The LQR design problem is presented in Section \ref{sec:problem}, along with its closed-form solution.
It is shown that the computation of the LQR state feedback gain by the closed-form solution requires knowledge of two  
large matrices: an extended observability matrix and a Toeplitz matrix of the plant's Markov parameters.
Then, in the ensuing sections \ref{sec:markov_estim} and \ref{sec:observ}, we provide algorithms to estimate these two matrices directly from data collected
from the plant. In Section \ref{sec:imc} we briefly review the Internal Model principle and the formulation of  
reference tracking as a state feedback problem. 
Two simulation examples are given in Section \ref{sec:simul} to illustrate the method's properties. 
One of our motivating applications -- the control of uninterruptible power sources --  is explored in Section  \ref{sec:ups},
where we present a 
practical application of our design methodology. It is seen
in the experimental results that our design compares
favorably with previously presented model-based solutions to this same practical problem.
Finally, concluding remarks are given in Section \ref{sec:conclusion}.

\section{PROBLEM STATEMENT} \label{sec:problem}
Consider a linear time-invariant discrete-time system
\begin{equation}\label{eq:sys}
\begin{aligned}
x(k+1) &= Ax(k)+Bu(k) \\
y(k) &= Cx(k)
\end{aligned}
\end{equation}
where $x(k)$ is an $n$-dimension state vector, $u(k)$ is a $p$-dimension input vector and $y(k)$ a $q$-dimension output vector.
The \textit{infinite horizon} LQR control problem can be summarized as follows:
find the optimal state feedback gain $K$ of the control law
\begin{equation}
u(k)=-Kx(k)
\label{eq:u(k)}
\end{equation}
such that the quadratic cost function
\begin{equation}
J = \sum_{k=0}^{\infty}{(y(k)^TQy(k)+u(k)^TRu(k))}
\label{eq:J_LQR}
\end{equation}
is minimized subject to system \eqref{eq:sys}, where $Q$ and $R$ are positive definite symmetric weighting matrices. The optimal gain is given by
\begin{equation}
K = (R + B^T P B)^{-1}(B^T P A)
\label{eq:opt_k}
\end{equation}
where $P$ is the unique positive definite solution to the discrete time algebraic Riccati equation (DARE)
\begin{multline}
P = A^T P A - (A^T P B)(R + B^T P B)^{-1} (B^T P A) \\ + C^TQC.
\label{eq:riccati}
\end{multline}
A closed-form solution to the DARE \eqref{eq:riccati} has been reported in \cite{lewis1981,furuta1993}. 
For a sufficient large $N$, this solution can be written as
\begin{equation}
P = \mathbf{O}^T(\mathbf{Q}_{N+1}^{-1}+\mathbf{S}\mathbf{R}_{N+1}^{-1}\mathbf{S}^T)^{-1}\mathbf{O},
\label{eq:cform}
\end{equation}
where
\begin{align}\label{eq:O}
\mathbf{O} =& \begin{bmatrix}C^T & (CA)^T & (CA^2)^T & \cdots & (CA^N)^T\end{bmatrix}^T,\\
\nonumber \mathbf{S} =& \begin{bmatrix}0 & \cdots & \cdots & \cdots & 0\\ 
CB & 0  & \ddots & & \vdots\\ CAB & CB & \ddots & \ddots & \vdots \\ 
\vdots & \vdots & \ddots & \ddots & \vdots\\ 
CA^{N-1}B & CA^{N-2}B & \cdots & CB & 0\end{bmatrix},
\end{align}
$$\mathbf{R}_j\triangleq\text{diag}(R,R,\ldots,R), ~~ \mathbf{Q}_j\triangleq\text{diag}(Q,Q,\ldots,Q),$$
with $\mathbf{R}_j$ and $\mathbf{Q}_j$ containing $j=N+1$ diagonal blocks each. 


The  matrix $\mathbf{O}$ is an \textit{extended observability matrix} for system \eqref{eq:sys} and $\mathbf{S}$ is a Toeplitz matrix of its Markov parameters:
\begin{equation}
M^{(i)} = CA^{(i-1)}B, ~~ i=1,2,\ldots,N .
\label{eq:markov_p}
\end{equation}
As shown in \cite{lewis1981}, using \eqref{eq:cform} in \eqref{eq:opt_k} and rearranging some terms, the 
LQR state feedback gain $K$ can be computed as a function of the Markov parameters as %
\begin{multline}
K = [R+\mathbf{M}^T(\mathbf{Q}_N^{-1}+\mathbf{S}\mathbf{R}_N^{-1}\mathbf{S}^T)^{-1}\mathbf{M}]^{-1}\\
\times \mathbf{M}^T(\mathbf{Q}_N^{-1}+\mathbf{S}\mathbf{R}_N^{-1}\mathbf{S}^T)^{-1}\mathbf{O}^{+},
\label{eq:K_dd}
\end{multline}
where
\begin{align}\label{eq:big_data}
\mathbf{M} &= \begin{bmatrix}M^{(1)} \\ M^{(2)}\\ \vdots \\ M^{(N)}\end{bmatrix} \in \Re^{qN\times p},\\
\mathbf{S} &= \begin{bmatrix}0 & \cdots & \cdots & \cdots & 0\\ 
M^{(1)} & 0  & \ddots & & \vdots\\ M^{(2)} & M^{(1)} & \ddots & \ddots & \vdots \\ 
\vdots & \vdots & \ddots & \ddots & \vdots\\ 
M^{(N-1)} & M^{(N-2)} & \cdots & M^{(1)} & 0\end{bmatrix} \in\Re^{qN\times pN}, \label{eq:S} \\
\mathbf{O}^{+} &= \begin{bmatrix} (CA)^T & (CA^2)^T & \cdots & (CA^N)^T\end{bmatrix}^T \in\Re^{qN\times n}.
\end{align}

Notice that gain $K$ in \eqref{eq:opt_k} is a function of $A, B,$ and $P$, which is also a function of system matrices $A,B$ and $C$. Now we have an expression that depends basically on $\mathbf{S}$ and $\mathbf{O}^+$. If these quantities could be obtained from data, then a data-driven method can be formulated. So, in order to succeed in this data-driven approach, we need to identify the system's Markov parameters -- the matrix $\mathbf{S}$ -- and an extended observability matrix $\mathbf{O}^+$.
%
%
%

%
%
%
Thus, let us pose the problem formally:

{\em
Given the data set
\begin{multline}
Z^T=[u(0), u(1),\ldots, u(T), y(0), y(1),\ldots, y(T),\\ x(0), x(1),\ldots, x(T)]
\label{eq:data_set}
\end{multline}
find the optimal state feedback gain $K$ as in \eqref{eq:K_dd}. To do so, a sequence of $N$ Markov parameters and the one-step ahead $N$ extended observability matrix must be estimated from data.
}

In the sequel we present a procedure to obtain both the Markov parameters and the extended observability matrix without using a model for the system. 

\section{Markov parameters estimation}
\label{sec:markov_estim}
We review next the estimation of the system Markov parameters via the so-called ARMarkov/Toeplitz models\footnote{We remark that estimating the system's Markov parameters is equivalent to identifying an $N$-th order FIR representation for the system.}. We follow the description in \cite{demoor:1988}. 


Let the states be repeatedly substituted $N+1$ times  in \eqref{eq:sys}; then
\begin{equation}\label{eq:sys_m}
\begin{aligned}
x(k+N+1) &= A^{N+1}x(k)+\mathbf{C}\mathbf{u}_m(k) \\
\begin{bmatrix} y(k) \\ y(k+1) \\ \vdots \\ y(k+N)\end{bmatrix} &= \mathbf{O}x(k)+\mathbf{S}\underbrace{\begin{bmatrix} u(k) \\ u(k+1) \\ \vdots \\ u(k+N)\end{bmatrix}}_{\mathbf{u}_m(k)}
\end{aligned}
\end{equation}
where $\mathbf{C}=[A^NB~\ldots~AB~B]$.

According to \cite{armakov}, as long as $(N+1)q\geq n$, it is guaranteed for an observable system that there exists a matrix $\mathbf{F}$ such that $A^{N+1}+\mathbf{F}\mathbf{O}=0$, which ensures that there exists an expression where the state is eliminated from \eqref{eq:sys_m}. This allows to write a predictor for the system's output as follows.

Let the Hankel matrix of a signal $e(k)$ be defined as
\begin{equation}
\resizebox{\hsize}{!}{$H(e(k)) \triangleq \begin{bmatrix}
e(k) & e(k+1) &  \cdots & e(k+L-1) \\
e(k+1) & e(k+2) & \cdots & e(k+L) \\
\vdots & \vdots & \ddots & \vdots \\
e(k+N-1) & e(k+N) & \cdots & e(k+N+L-2) 
\end{bmatrix}$.}
\label{eq:hankel}
\end{equation}

Define the set of data matrices 
\begin{equation}\label{eq:data_matrices}
\begin{aligned}
U_p=H(u(0)) & ~~~~~~~~~ U_f=H(u(N)) \\
Y_p=H(y(0)) & ~~~~~~~~~ Y_f=H(y(N)).
\end{aligned}
\end{equation}
Then a predictor of the system output can be written as \cite{demoor:1988}
\begin{equation}\label{eq:sys_m_elimx}
Y_f = \underbrace{[\mathbf{O}(\mathbf{C}+\mathbf{F}\mathbf{S})~~-\mathbf{O}\mathbf{F}~~\mathbf{S}]}_{\mathbf{W}}\underbrace{\begin{bmatrix}U_p \\ Y_p \\ U_f\end{bmatrix}}_{\mathbf{\Phi}}.
\end{equation}


Thus, an estimate $\widehat{\mathbf{W}}$ of $\mathbf{W}$ can be obtained by solving the least-squares problem
\begin{equation}
\mathbf{W}=Y_f\mathbf{\Phi}^{\dagger}
\label{eq:MQ}
\end{equation}
with $[\cdot]^{\dagger}$ denoting the Moore-Penrose pseudo-inverse, and an estimate $\widehat{\mathbf{S}}$ as the rightmost $pN$ columns of $\widehat{\mathbf{W}}$.
%
%
%
%
%
Moreover, $\mathbf{M}$ can be extracted from the first column of $\widehat{\mathbf{S}}$.
This estimation has been shown to be consistent for the Markov parameters \cite{kamrunnahar2000}.

For stable systems, the choice of $N$ is closely related to the system's open loop settling time, as the contribution of $A^N$ decreases as $N$ increases. Also, $L$ must satisfy $L\geq 3q(N+1)$ so \eqref{eq:MQ} has a solution. 

\section{Extended observability matrix estimation} \label{sec:observ}
Since a state feedback control is to be implemented, we can assume that the states are measurable. Hence, in this section we present two original algorithms to identify an extended observability matrix in the same state coordinates we are measuring and later we discuss their properties. We define the vector of measured states by 
\begin{equation}
\mathbf{X}\triangleq [x(0)~~x(1)~~x(2)~~\cdots~~x(L-1)].
\label{eq:x_hankel}
\end{equation}

\subsection{Algorithm 1}




The output equation \eqref{eq:sys_m} can also be written with Hankel matrices without eliminating the state vector. It can be put in an extended output matrix equation  as \cite{demoor:1988}
\begin{equation}
Y_p = \mathbf{O}\mathbf{X}+\mathbf{S}U_p.
\label{eq:sys_hankel}
\end{equation}

Since an estimate $\widehat{\mathbf{S}}$ can be obtained using the algorithm provided in Section \ref{sec:markov_estim}, and $U_p$ and $Y_p$ can be formed with collected data, we can then solve the system of equations in \eqref{eq:sys_hankel} for $\mathbf{O}$:
\begin{equation}
\widehat{\mathbf{O}} = (Y_p-\widehat{\mathbf{S}}U_p)\mathbf{X}^{\dagger}.
\label{eq:observ}
\end{equation}
To obtain $\mathbf{O}^{+}$, one simply removes the first $q$ rows of $\mathbf{O}$. 

%
%
%
\subsection{Algorithm 2}
First, define $U_{po}$ as the geometric operator that projects the row space of a matrix onto the orthogonal
complement of the row space of the matrix $U_{p}$
\begin{equation}
U_{po} \triangleq I_{L} - U_p^T(U_pU_p^T)^{-1}U_p
\label{eq:ortho_u}
\end{equation}
where $I_L$ is an identity matrix of size $L$. Then by post-multiplying the extended output matrix \eqref{eq:sys_hankel} by $U_{po}$, we have
\begin{align} \nonumber
Y_p U_{po} &= \mathbf{O}\mathbf{X}U_{po}+\mathbf{S}U_p U_{po} \\
Y_p U_{po} &= \mathbf{O}\mathbf{X}U_{po}
\label{eq:sys_hankel_ortho}
\end{align}
Notice that by using the projection $U_{po}$ we eliminate the need to know, or estimate, $\mathbf{S}$.

An estimator of the extended observability matrix can thus be computed as
\begin{equation}
\widehat{\mathbf{O}} = (Y_p U_{po})(\mathbf{X}U_{po})^{\dagger}.
\label{eq:observ2}
\end{equation}

\subsection{Estimates properties}

The estimates just provided allow the computation of the state feedback gain according to \eqref{eq:K_dd} and, as shown in \cite{lewis1981}, 
the gain thus computed converges asymptotically, as $N\rightarrow\infty$, 
to the optimal LQR state feedback. A simulation example in Section \ref{sec:simul} illustrates this
property.

On the other hand, when the state measurement is corrupted by noise, one can expect some bias in the estimates
of the extended observability matrix, since both solutions presented are of a least squares nature. We now briefly
discuss the bias and covariance of these estimates. Consider the system state-space representation \eqref{eq:sys}
with the noise terms
\begin{equation}\label{eq:sys_noise}
\begin{aligned}
x(k+1) &= Ax(k)+Bu(k)+Ev(k) \\
y(k) &= Cx(k)+Fw(k)
\end{aligned}
\end{equation}
with $v(k)$ and $w(k)$ white noise sequences. We can write the extended output equation as
\begin{equation}
Y_p = \mathbf{O} \mathbf{X} + \mathbf{S} U_p + \mathbf{S}_E V_p + F W_p
\label{eq:ext_out_noise}
\end{equation}
where $V_p=H(v(0))$ and $W_p=H(w(0))$ are Hankel matrices of the noise sequences $v(k)$ and $w(k)$ respectively, and $\mathbf{S}_E$ has the same structure as $\mathbf{S}$ with $E$ 
\textit{in lieu} of $B$.
We can also write
\begin{equation}
\mathbf{X} = \mathbf{X}^o+E\mathbf{V}
\label{eq:x_noise}
\end{equation}
where $\mathbf{V}$ is a row vector of the noise sequences $v(k)$. We also assume the sequence $u(k)$ uncorrelated with $v(k)$.

\subsubsection{Algorithm 1}
Let $\Ex[\cdot]$ denote the expected value function. The bias of the first algorithm is given by
\begin{equation}
\mathcal{B}(\widehat{\mathbf{O}})=\Ex[\widehat{\mathbf{O}}]-\mathbf{O} = \Ex[(Y_p-\widehat{\mathbf{S}}U_p)\mathbf{X}^{\dagger}-\mathbf{O}]. \label{aux1}
\end{equation}

Inserting \eqref{eq:ext_out_noise} into \eqref{aux1}
\begin{align} \nonumber
\mathcal{B}(\widehat{\mathbf{O}})={}& \Ex[(\mathbf{O} \mathbf{X} + \mathbf{S} U_p + \mathbf{S}_E V_p + F W_p-\widehat{\mathbf{S}}U_p)\mathbf{X}^{\dagger}-\mathbf{O}] \\
\begin{split}={}&\Ex[\mathbf{O} (\mathbf{X}\mathbf{X}^{\dagger}-I) + (\mathbf{S}-\widehat{\mathbf{S}})U_p\mathbf{X}^{\dagger} + \mathbf{S}_E V_p \mathbf{X}^{\dagger} \\ &  + F W_p \mathbf{X}^{\dagger}] \end{split} \label{aux2}
\end{align}

First term of \eqref{aux2} is null since $\mathbf{X}$ has real entries, and the second term is null either if $\widehat{\mathbf{S}}=\mathbf{S}$ or because $u(k)$ is uncorrelated with $v(k)$. We then have
\begin{equation}
\mathcal{B}(\widehat{\mathbf{O}})=\mathbf{S}_E\Ex[V_p \mathbf{X}^{\dagger}] + F \Ex[W_p \mathbf{X}^{\dagger}].
\label{eq:algo1b1}
\end{equation}
If we further assume that $v(k)$ is uncorrelated with $w(k)$, then
\begin{equation}
\mathcal{B}(\widehat{\mathbf{O}})=\mathbf{S}_E\Ex[V_p \mathbf{X}^{\dagger}].
\label{eq:algo1b2}
\end{equation}

The covariance of the first algorithm is given by
\begin{align}
\mathcal{V}(\widehat{\mathbf{O}})={}&\Ex[(\widehat{\mathbf{O}}-\Ex[\widehat{\mathbf{O}}])(\widehat{\mathbf{O}}-\Ex[\widehat{\mathbf{O}}])^T] \\ \nonumber
\begin{split}={}&\Ex[((\mathbf{S}-\widehat{\mathbf{S}})U_p\mathbf{X}^{\dagger} + \mathbf{S}_E V_p \mathbf{X}^{\dagger} + F W_p \mathbf{X}^{\dagger}
\\ & -  \mathbf{S}_E \Ex[V_p \mathbf{X}^{\dagger}] - F\Ex[ W_p \mathbf{X}^{\dagger}])
\\ & \times ((\mathbf{S}-\widehat{\mathbf{S}})U_p\mathbf{X}^{\dagger} + \mathbf{S}_E V_p \mathbf{X}^{\dagger} + F W_p \mathbf{X}^{\dagger}
\\ & -  \mathbf{S}_E\Ex[ V_p \mathbf{X}^{\dagger}] - F\Ex[ W_p \mathbf{X}^{\dagger}])^T] 
\end{split} 
\end{align}
which after some algebraic manipulation results in
\begin{align}
\begin{split}\label{eq:algo1var}
\mathcal{V}(\widehat{\mathbf{O}})={}&\Ex\left\{[(\mathbf{S}-\widehat{\mathbf{S}})U_p\mathbf{X}^{\dagger} + \mathbf{S}_E V_p \mathbf{X}^{\dagger} + F W_p \mathbf{X}^{\dagger}]\Lambda_1^T\right.
\\ & + \left.[(\mathbf{S}-\widehat{\mathbf{S}})U_p\mathbf{X}^{\dagger}+\Lambda_{1}][(\mathbf{S}-\widehat{\mathbf{S}})U_p\mathbf{X}^{\dagger}]^T\right\}
\end{split} \\
\Lambda_1={}&\mathbf{S}_E (V_p \mathbf{X}^{\dagger}-\Ex[V_p \mathbf{X}^{\dagger}])+F (W_p \mathbf{X}^{\dagger}-\Ex[W_p \mathbf{X}^{\dagger}])
\end{align}



\subsubsection{Algorithm 2}

Following the same steps as in Algorithm 1, the following expressions are obtained for the bias and covariance of the
estimates given by Algorithm 2.  
\begin{equation}
\mathcal{B}(\widehat{\mathbf{O}})=\mathbf{S}_E\Ex[(V_p U_{po})(\mathbf{X}U_{po})^{\dagger}].
\label{eq:algo2b2}
\end{equation}

\begin{align}
\begin{split}\label{eq:algo2var}
\mathcal{V}(\widehat{\mathbf{O}})={}&\Ex\left\{\mathbf{B}[\mathbf{B} + \mathbf{S}_E V_p \mathbf{X}^{\dagger} + F W_p \mathbf{X}^{\dagger}]^T\right. \\
+ &\left.[\mathbf{S}_E V_p \mathbf{X}^{\dagger} + F W_p \mathbf{X}^{\dagger}]\mathbf{B}^T\right\}
\end{split} 
\end{align}
\begin{multline}
\Lambda_{2}=\mathbf{S}_E[(V_p U_{po})(\mathbf{X}U_{po})^{\dagger}-\Ex[(V_p U_{po})(\mathbf{X}U_{po})^{\dagger}]]
\\+  F[(W_p U_{po})(\mathbf{X}U_{po})^{\dagger}-\Ex[(W_p U_{po})(\mathbf{X}U_{po})^{\dagger}]]
(\mathbf{X}U_{po})^{\dagger}
\end{multline}

Notice that, as expected, the bias of both estimates is inversely proportional to the signal to noise ratio and the estimates will
be unbiased only if there is no noise in the measurement. We provide an illustrative example in Section \ref{ssec:ex_noise} to compare the bias, covariance and -- most importantly -- the mean square error resulting from the two algorithms.


%
%
%

\section{Internal model controller and augmented state space} \label{sec:imc}
Most practical control applications consider the reference tracking problem. Thus, in this section we show how to use open-loop data to obtain an augmented state and output vectors in order to adjust the gains also for reference tracking considering a feedback loop with an internal model controller (IMC) \cite{francis1976}.
Let 
\begin{equation}
\mathcal{C}=(A_c,B_c)
\label{eq:imc}
\end{equation}
be a discrete-time realization of the internal model controller state equation, whose $n_c$ states are measurable. Assume that every output of the system is to follow a reference represented by the internal model controller. Then the augmented open-loop space-state representation of system \eqref{eq:sys} with controller \eqref{eq:imc} is given by
\begin{equation}
\begin{aligned}
x_a(k+1) &= \underbrace{\begin{bmatrix}A & \mathbf{0}_{n,n_c\times q}\\ -C\otimes B_c & I_q\otimes A_c \end{bmatrix}}_{A_a}x_a(k)+\underbrace{\begin{bmatrix}B \\ \mathbf{0}_{n_c\times q,p} \end{bmatrix}}_{B_a}u(k) \\
y_a(k) &= \underbrace{\begin{bmatrix}C & \mathbf{0}_{q,n_c\times q}\\ \mathbf{0}_{n_c\times q,n} & I_{n_c\times q}\end{bmatrix}}_{C_a}x_a(k)
\end{aligned}
\label{eq:aug_sys}
\end{equation}
where $\otimes$ is the Kronecker product. 

From partition $(2,1)$ of $A_a$ we see that the open-loop IMC states needed to compute the gain $K$ can be obtained by simply filtering the plant outputs by $-\mathcal{C}$.

For example, the \textit{state equation} of the integrator $\frac{1}{z-1}$ can be represented by
\begin{equation}
x_c(k+1) = \underbrace{1}_{A_c} x_c(k) + \underbrace{1}_{B_c} (r(k)-y(k))
\label{eq:integral}
\end{equation}
and a resonant controller at frequency $\omega_n$ with a pre-warping Tustin representation can be realized as
\begin{equation}
x_c(k+1)\underbrace{\begin{bmatrix}0 & 1 \\ -1 & 2\cos(\omega_nT_s)\end{bmatrix}}_{A_c} x_c(k)+\underbrace{\begin{bmatrix}0 \\ 1 \end{bmatrix}}_{B_c}(r(k)-y(k))\\
\label{eq:resonant}
\end{equation}
where $T_s$ is the sampling time.

Let $x_{IMC}(t)$ represent the state of the open-loop IMC\footnote{which can be obtained with MATLAB command \texttt{lsim}; for instance, \texttt{xIMC=-lsim(IMC,y)}.}, and the augmented output and state vector be given by $y_a(t) = [y(t)~ x_{IMC}(t)]^T$ and $x_a(t) = [x(t) ~x_{IMC}(t)]^T$, respectively. These are the vectors that should be used, along with $u(t)$, to estimate the Markov parameters and the extended observability matrix.

\section{Simulation Examples}\label{sec:simul}

\subsection{Regulation control}
In this first example we illustrate the convergence of the proposed method.
Consider a system as in \eqref{eq:sys} whose matrices are given by
\begin{equation}\label{eq:ex_sys}
A = \begin{bmatrix} 1 & 0.15 \\ -0.2 & 0.6 \end{bmatrix}~~B = \begin{bmatrix} 0.04 & 0.01 \\ 0.02 & -0.01 \end{bmatrix}~~ C= \begin{bmatrix} 1 & 2 \\ 0 & 1 \end{bmatrix}\cdot
\end{equation}
The sampling time is $T_s=1$~s. Also, let the performance requirements be given by
$R = 0.2I_q$ and $Q = 20I_q$, where $I_q$ is the identity matrix of size $q$, that is, we are valuing more the evolution of the states than the control effort.  The \textit{model-based} optimal LQR controller is given by
\begin{equation}
K=\begin{bmatrix} 4.6491 & 7.5226\\ 1.4461 & -1.9886 \end{bmatrix}\cdot
\label{eq:K_modelbased}
\end{equation}

In order to apply the proposed methodology we set an open-loop experiment where the input signal is a PRBS with amplitude $1$ and length $1022$ and the output data was collected, i.e., both $y(k)$ and $x(k)$. We identified two state-feedback gains: one with $N=10$ and the other with $N=50$, the latter been close to the system's open loop settling time (approximately $45$ samples). We obtained $\hat{K}_{10}=\begin{bmatrix} 4.2314 & 7.644 \\ 1.127 & -1.8959 \end{bmatrix}$ and 
$\hat{K}_{50}=\begin{bmatrix} 4.6491 & 7.5226\\ 1.4461 & -1.9886 \end{bmatrix}$.


Notice that with $N=50$, $\hat{K}$ equals the model-based solution \eqref{eq:K_modelbased} up to the fifth significant digit. 
\subsection{The noisy case: observability matrix properties}
\label{ssec:ex_noise}
We provide now a simple example to illustrate bias and covariance of the observability matrix estimators. Consider a system with state-space matrices
\begin{equation}
A = 0.14, ~~~ B = 1.72, ~~~ C=1, ~~~ E=1, ~~~ F=0,
\label{eq:ex_sys_noise}
\end{equation}
and let the LQR performance matrices be $R = 0$ and $Q = 1$, so we are aiming for a dead-beat control, and convergence can be found with $N_m=3$ (due to row removal, that means the extended observability matrix will be of size $2$). The actual extended observability matrix is then $\mathbf{O} = [0.14~~0.0196]^T$.

\begin{figure}[htb]%
\centering
\includegraphics[width=0.9\columnwidth]{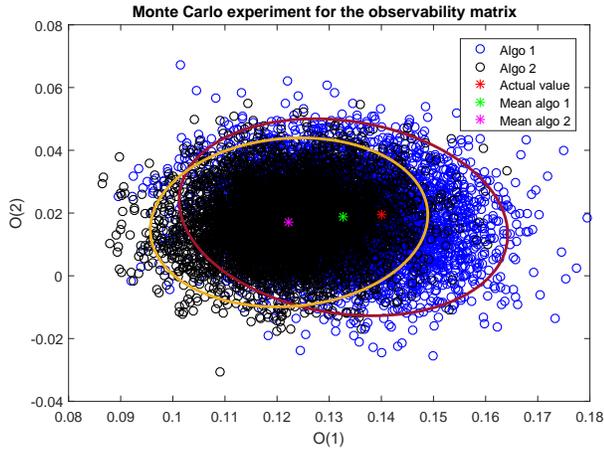}%
\caption{Estimates of the observability matrix in a Monte Carlo experiment. The ellipses represent the covariance regions around the mean value with $95\%$ confidence.}%
\label{fig:observ_prop}%
\end{figure}

We set a Monte Carlo experiment with $5000$ runs and with a PRBS input of length $1022$ and $e(t)$ as white noise with variance $\sigma^2=0.1$. Fig. \ref{fig:observ_prop} portrays the results obtained with both estimators \eqref{eq:observ} and \eqref{eq:observ2}.

As mentioned before, the estimate bias with Algorithm $1$ is smaller, whereas with Algorithm $2$ a smaller covariance is achieved. In fact, we obtained $\Ex(\widehat{\mathbf{O}}_1) = [0.1327~~0.01865]^T$ and $\Ex(\widehat{\mathbf{O}}_2) = [0.1223~~0.01704]^T$, and the eigenvalues $\lambda$ of the covariance matrices $\lambda[\Ex(\widehat{\mathbf{V}}_1)] = [1.365~~1.938]^T\times 10^{-4}$ and $\lambda[\Ex(\widehat{\mathbf{V}}_2)] = [1.093~~1.3]^T\times 10^{-4}$. Notice that the largest eigenvalue with Algorithm $2$ is approximately the smallest eigenvalue with Algorithm $1$.

We also computed the eigenvalues of the MSE matrix of both algorithm and obtained $\lambda(MSE_1)= [1.723\times 10^{-4}~~1.812\times 10^{-2}]^T$ and $\lambda(MSE_2)= [1.182\times 10^{-4}~~1.537\times 10^{-2}]^T$. Note that Algorithm $2$ provides much smaller MSE.

\section{Control of an UPS} \label{sec:ups}
We now consider a practical application of the proposed methodology to an uninterruptible power supply (UPS).  This plant has been studied before in \cite{pereira2014, lorenzini2015}.

\begin{figure}[htb]%
\centering
\includegraphics[width=0.8\columnwidth]{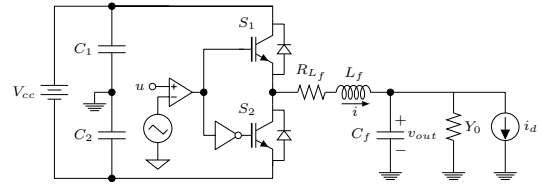}%
\caption{Schematic representation of the UPS with load.}%
\label{fig:ups_scheme}%
\end{figure}

Consider the simplified electrical diagram of the output
stage of a single-phase UPS system, as illustrated in Figure \ref{fig:ups_scheme}. 

\setlength{\arraycolsep}{3pt}

The load effect on the system output is modeled by a parallel connection of an uncertain admittance $Y_0(t)$ and an unknown periodic disturbance given by the current source $i_d(t)$. The PWM (Pulse Width Modulation) comparator input is modeled by a gain $K_{PWM}$ multiplied by the control input. Also, defining the system state vector as the inductor current and the capacitor voltage, $x(t) = [i(t)~~v(t)]^T$, the \textit{continuous-time} state-space representation for the UPS system is given by:
\begin{equation}\label{eq:ex_ups}
\begin{aligned}
\dot{x}(t) &= \underbrace{\begin{bmatrix} \frac{-R_{L_f}}{L_f} & \frac{-1}{L_f} \\[4pt] \frac{1}{C_f} & \frac{-Y_0(t)}{C_f} \end{bmatrix}}_{A}x(t)+\underbrace{\begin{bmatrix} \frac{K_{PWM}}{L_f} \\[4pt] 0 \end{bmatrix}}_{B}u(t)+ \underbrace{\begin{bmatrix} 0 \\[4pt] \frac{-1}{C_f} \end{bmatrix}}_{B_d}i_d(t) \\
y(t) &= \underbrace{[~~~ 0 ~~~~~~~ 1 ~~~]}_{C} x(t)
\end{aligned}
\end{equation}
where $u(t)$ is the PWM control input
and $y(t)$ is the output voltage to be controlled.

\setlength{\arraycolsep}{6pt}
Closed-loop reference is typically a sinusoid and for our case we have $r(t)=127\sqrt{2}\sin(120\pi t)$.
Since the reference signal is a sinusoid, then the right choice of the IMC is a resonant controller \eqref{eq:resonant}.
Admittance $Y_0(t)$ can be set as an open circuit (no load), as a nominal resistance $R_0 = 6~\Omega$ 
and as a non-linear load given by a full-bridge circuit.

%

The control objective can be summarized as: \textit{design a data-driven state feedback controller for sinusoid reference tracking for the UPS operating with \textbf{non-linear load}}.

\begin{figure}[htb]%
\includegraphics[width=0.98\columnwidth]{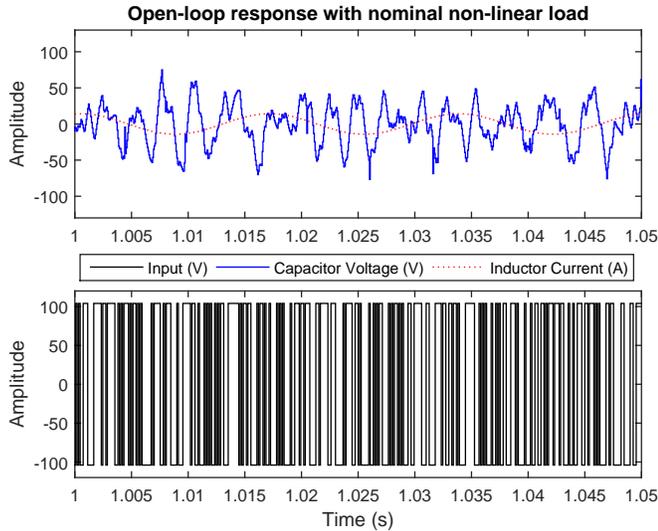}%
\caption{Open-loop response with nominal non-linear load.}%
\label{fig:ups_ma}%
\end{figure}

In order to obtain meaningful data from the system, we set an open-loop experiment as follows. The sampling time was set to $T_s=1/15000$~s; the input of the PWM was set as a PRBS with amplitude $\pm 104$~V and with length $75000$ samples (i.e., a $5$ seconds signal); current and voltage were measured as portrayed (zoomed time scale) in Fig. \ref{fig:ups_ma} for the UPS operating with its \textit{non-linear load}. The output voltage was filtered by \eqref{eq:resonant} to obtain the IMC states for our proposed algorithm.

The parameters chosen for the LQR were $Q=200I_3$ and $R=5000$, that is, we strongly penalized the control signal as to try to achieve closed-loop stability even when there is no load in the UPS and to reduce sensibility due to noise, specially for the current measurement. Prior to any knowledge about the system settling time, we also selected $N=150$. The obtained LQR gain is
\begin{equation}
K = [4.85548~~5.54514~~0.638479~~-0.644019].
\label{eq:K_ups}
\end{equation}

%
%

\begin{figure}[htb]%
\centering
\includegraphics[width=0.98\columnwidth]{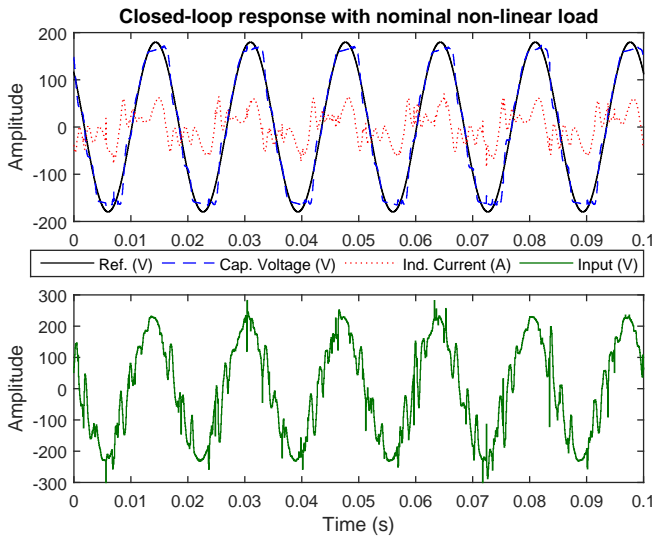}%
\caption{Closed-loop response with nominal non-linear load.}%
\label{fig:ups_nominal}%
\end{figure}

Fig. \ref{fig:ups_nominal} shows the closed-loop response for the UPS operating at nominal capacity with non-linear load.
Stability and reference tracking were achieved with corresponding Total Harmonic Distortion (THD) of $11.7\%$ -- result similar to the one obtained in \cite{pereira2014}, in which the state feedback gain was designed using a full plant model and a Linear Matrix Inequality approach.

\begin{figure}[htb]
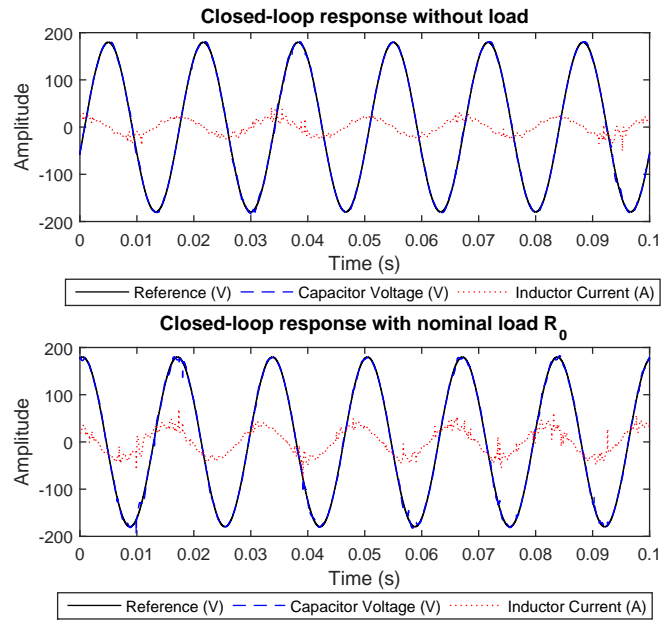
%
\centering
\includegraphics[width=0.98\columnwidth]{ups_mf_vazio.eps}%
\vspace{2pt}
\includegraphics[width=0.98\columnwidth]{ups_mf_lin.eps}%
\caption{Closed-loop response with no output load and with nominal linear load.}%
\label{fig:ups_vazio}%
\end{figure}

We also applied the obtained controller to different scenarios: $(a)$ the open circuit (no load) case and $(b)$ the nominal linear load $R_0$. Fig. \ref{fig:ups_vazio} 
shows the closed-loop responses for cases $(a)$ and $(b)$ respectively. Closed-loop stability was achieved and with very small THD -- $1.4\%$ and $1.8\%$ respectively --, even though the controller was not designed for these scenarios.




Notice that with this approach we obtained a linear state-feedback gain with only one experiment on the plant, even though the actual plant has a strong nonlinear behavior and a single linear model would not describe the system with reasonable accuracy. If data were used to identify a plant model, then more than one experiment would be necessary in order to evaluate a plant model and an uncertainty matrix.


\section{Conclusion} \label{sec:conclusion}
In this paper we provided a data-driven method to compute the infinite horizon LQR state feedback gain, without identifying
a model of the plant. In our method, the feedback gain is computed from a batch of data and converges to the infinite
horizon LQR gain as the amount of data and, by consequence, the number of estimated Markov parameters grow. Simulation examples illustrated the convergence of the method and an
experimental application to an UPS showed its practical applicability.

\bibliographystyle{IEEEtran}
\bibliography{baseref}

\end{document}